\documentclass[a4paper]{amsart}

\usepackage{amsmath}
\usepackage{amsfonts}
\usepackage{amssymb}
\usepackage{graphicx}
\newtheorem{theo}{Theorem}[section]
\theoremstyle{plain}

\newtheorem{cor}[theo]{Corollary}

\newtheorem{example}[theo]{Example}

\newtheorem{lemma}[theo]{Lemma}

\newtheorem{proposition}[theo]{Proposition}
\newtheorem{remark}[theo]{Remark}

\numberwithin{equation}{section}
\begin{document}
\baselineskip=17pt 
\title[separate and joint continuity]{Some results on separate and joint
continuity}
\author{A. Bareche}
\address{Universit\'e de Rouen, 
UMR CNRS 6085, Avenue de l'Universit\'e, BP.12, F76801 Saint-\'Etienne-du-Rouvray, France.}
\email{aicha.bareche@etu.univ-rouen.fr}
\author{A. Bouziad}
\address[]
{D\'epartement de Math\'ematiques,  Universit\'e de Rouen, UMR CNRS 6085, 
Avenue de l'Universit\'e, BP.12, F76801 Saint-\'Etienne-du-Rouvray, France.}
\email{ahmed.bouziad@univ-rouen.fr}
 

\begin{abstract}  Let $f: X\times K\to \mathbb R$ be a separately continuous
function and $\mathcal C$ a countable collection of  subsets of $K$. Following
a result of  
Calbrix and  Troallic, there is a residual set of points $x\in X$
such that $f$ is jointly continuous at each point of $\{x\}\times Q$, where $Q$
is the set of $y\in K$ for which  the collection
$\mathcal C$ includes a  basis of neighborhoods in $K$.  The
particular
case when the factor $K$ is second countable was recently extended by Moors and
Kenderov to any \v Cech-complete Lindel\"of space $K$ and Lindel\"of
$\alpha$-favorable
$X$, improving   a generalization
of Namioka's theorem obtained by  Talagrand. 
Moors proved the same result when   $K$  is  a Lindel\"of $p$-space  
and $X$ is conditionally $\sigma$-$\alpha$-favorable space. 
Here we add new  results of this sort  when the factor $X$
is  $\sigma_{C(X)}$-$\beta$-defavorable and when the assumption ``base
of
neighborhoods" 
in 
Calbrix-Troallic's result is replaced by a type of  countable
completeness. 
    The paper also provides further information
about the class of Namioka spaces.
\end{abstract} 
\maketitle
\vskip 4mm
\section{Introduction}

If $K$, $X$ are topological spaces, a mapping $f:X\times K\to \mathbb R$ is said
to be separately
continuous if for every $x\in X$ and $y\in K$,
the mappings $f(x,.): K\to\mathbb R$ and $f(.,y): X\to\mathbb R$ are continuous,
the reals being equipped
with the usual topology. 
The spaces $K$ and $X$ satisfy the Namioka property  ${\mathcal N}(X,K)$ 
if every separately continuous map $f:X\times K\to\mathbb  R$ is 
(jointly) continuous at each point of a subset of $X\times K$
of the form $R\times K$, where $R$ is a dense subset of $X$ \cite{nam}.
Following \cite{chr}, the space $X$ is called
a Namioka space if the property ${\mathcal N}(X,K)$
holds for every compact $K$. It is well known that every Tychonoff Namioka space
is a Baire space \cite{sr2}. 
Following \cite{de1}, a compact space $K$ is said to be
co-Namioka  if ${\mathcal N}(X,K)$ holds for every
Baire space $X$. The class of co-Namioka spaces contains several classes of
compact spaces appearing
in Banach spaces theory, like Eberlein or Corson compactums (\cite{dev},
\cite{de2}); in this connection, the reader is referred to
\cite{merneg,pio, bormoo} and the references therein for more
information.  On the other hand, every $\sigma$-$\beta$-defavorable space (see
below)
is a  Namioka space; this is Christensen-Saint Raymond's theorem
\cite{chr,sr2}. It is also well known that
within  the class
of metrizable or separable spaces, Namioka spaces and
$\sigma$-$\beta$-defavorable
spaces coincide \cite{sr2},  a result that we will   improve below by extending
it to 
  Grothendieck-Eberlein spaces (see also Proposition 5.5). Any
Baire 
space which is a $p$-space (in Arhangel'ski\v \i's sense) or
 $K$-analytic  
is $\sigma$-$\beta$-defavorable, hence a Namioka space; see respectively
\cite{bo2} and \cite{de1}. It should be noted  that
the method of \cite{de1}
can be used to extend this result of Debs to any Baire space which is  dominated
 by the irrationals in the sense of \cite{tk2}. In addition, a Baire space which is game
 determined in the sense of Kenderov and Moors in \cite{kenmor} is $\sigma$-$\beta$-defavorable. In particular, a Baire
 space which has countable separation is $\sigma$-$\beta$-defavorable. \par 
The class
of $\sigma$-$\beta$-defavorable spaces   
is defined in term of a topological game $\mathcal J$ introduced
(in a strong form) by Christensen \cite{chr} and later modified by Saint Raymond
in \cite{sr2}. In the game $\mathcal J$ on  the space $X$, two players $\alpha$ 
and $\beta$ choose alternatively
a decreasing sequence $V_0\supseteq U_0\supseteq \ldots\supseteq V_n\supseteq
U_n\ldots$ of nonempty 
open subsets of $X$ and a sequence $(a_n)_{n\in\mathbb N}\subset X$ as follows:
Player $\beta$ moves first 
and chooses
$V_0$; then Player $\alpha$ gives $U_0\subset V_0$
and $a_0\in X$. At the $(n+1)$th step, Player $\beta$
chooses an open set $V_{n+1}\subset U_n$ then Player $\alpha$ responds by giving
$U_{n+1}\subset V_{n+1}$ and $a_{n+1}\in X$. 
The {\it play} $(V_n, (U_n,a_n))_{n\in\mathbb N}$ is won by Player
$\alpha$ if   $$(\cap_{n\in \mathbb N}U_n)\cap \overline{\{a_n
:n\in\mathbb N\}}\not=\emptyset.$$  
The space $X$ is said to be $\sigma$-$\beta$-defavorable if there
is no winning strategy for Player $\beta$ in the game  $\mathcal J$.
 \par 
The problem of knowing to what extent can we weaken the assumption of
compactness on the factor $K$ has interested several authors. In this work, we
are interested in certain results obtained on this issue, that we describe now. 
Let $(U_n)_{n\in\mathbb N}$ be a sequence
of  open subsets of $K$. In \cite{caltro}, 
 Calbrix
and  Troallic have shown  that there is a residual
set $R\subset X$ such that the
separately continuous mapping
$f:X\times K\to\mathbb R$
is continuous at each point of $R\times Q$, where
$Q$ is the set of points $x\in K$ admitting a 
subsequence of $(U_n)_{n\in\mathbb N}$ as a neighborhoods basis. In particular,
the property ${\mathcal N}(X,K)$ holds for every
second countable space $K$ and every Baire space $X$. A similar result has been
proved previously by Saint Raymond \cite{sr1} in the case where $K$ and $X$ are
both Polish.  In the same direction,  Talagrand has demonstrated in \cite{tal} 
 that ${\mathcal N}(X,K)$
 holds when $K$ is \v Cech-complete Lindel\"of 
and $X$ is compact (also announcing  the same result for $X$  \v Cech-complete
complete).
Mercourakis     and Negrepontis suspected in their  article \cite{merneg} the
possibility of extending these results  in case where $K$ is 
Lindel\"of $p$-space, which has been established with success by Moors in a
recent article \cite{moo}
assuming $X$ to be ``conditionally'' $\alpha$-favorable. Shortly before that,
Moors and Kenderov extended in \cite{kenmoo} Talagrand's theorem    
to every $\alpha$-favorable Lindel\"of space $X$. 
As the class of $\sigma$-$\beta$-defavorable spaces encompasses so nicely
different
types of Namioka spaces, it seemed to us that it would be interesting to know if
 some results of this kind remain valid in the framework of this class.\par
The basic idea here is the reuse of the approach in \cite{bo1},  where
a simplified proof is given for Christensen-Saint Raymond's theorem. In Theorem
3.2,
the result of Calbrix and Troallic is considered in a more general
configuration, 
replacing the set $Q$ by the set of $x\in K$
for which there is  a subsequence of $(U_n)_{n\in\mathbb N}$
containing $x$ and satisfying a sort of 
countable completeness (the precise definition is given in Section 3). This also
 concerns the
above result by Moors. In Theorem 3.1, we shall examine the case where the
sequence $(U_n)_{n\in\mathbb N}$ is  a sequence of countable 
(not necessarily open) covers of $K$, which will allow us to unify the result of
Talagrand (including the $K$-analytic variant of his theorem)  and that of
Kenderov and Moors.  
 Concerning the factor $X$, we shall do a functional  
adjustment to the game of Christensen-Saint Raymond,  thereby obtaining a class
wider than that of $\sigma$-$\beta$-defavorable spaces  whose members are still
Namioka spaces. For instance, this new class contains
all pseudocompact spaces. Related to this last result, a more general statement
is proved
in Proposition 5.5 in Section 5 which includes  some additional results and
observations.
\section{Functional variants of  Christensen-Saint Raymond's game} 
\vskip 2mm
        \noindent {\bf  The game ${\mathcal  J}_\Gamma$:} Let $\Gamma\subset
{\mathbb R}^X$. The
         game ${\mathcal J}_\Gamma$ differs from the game  $\mathcal J$ only
        in the winning condition: Player
$\alpha$ is declared to be the winner of the play  $(V_n,(U_n,a_n))_{n\in\mathbb
N}$ if for each $g\in
\Gamma$ there exists  $t\in \cap_{n\in \mathbb N}U_n$ 
such that  $$g(t)\in\overline{\{g(a_n): n\in\mathbb N\}}.$$
The space $X$ is said to be  {\it $\sigma_\Gamma$-$\beta$-defavorable} if Player
$\beta$ has no winning strategy in the game ${\mathcal  J}_\Gamma$. Using a
terminology from \cite{moo}, we shall
say that
 $X$ is  {\it conditionally $\sigma_\Gamma$-$\alpha$-favorable} if Player
$\alpha$
has a strategy  $\tau$ so that for any compatible play
$(V_n,(U_n,a_n))_{n\in\mathbb N}$ 
satisfying  $\cap_{n\in\mathbb N}U_n\not=\emptyset$,   for every
$g\in \Gamma$ there is  $t\in\cap_{n\in\mathbb N}U_n$ such that
 $$g(t)\in\overline{\{g(a_n):
n\in\mathbb N\}}.$$

It will be   useful for our purpose to consider the closely related game
${\mathcal J}_{\Gamma}^*$ where   Player $\alpha$ has not to produce
the sequence $(a_n)_{n\in\mathbb N}\subset X$ but wins the play
$((V_n,U_n))_{n\in\mathbb N}$ if (and only if)
for each sequence $(a_n)_{n\in\mathbb N}$ such that
 $a_n\in U_n$ ($n\in\mathbb N$) and for each $g\in \Gamma$, there exists $t\in
\cap_{n\in\mathbb N} U_n$ such that
$$g(t)\in\overline{\{ g(a_n):n\in\mathbb N\}}.$$  
We make similar definitions with ${\mathcal J}_{\Gamma}$ replaced by ${\mathcal
J}_{\Gamma}^*$;
for instance, $X$ is said to be  $\sigma_{\Gamma}^*$-$\beta$-defavorable space
if Player $\beta$ has no winning strategy in the game ${\mathcal J}_\Gamma^*$.
  \par
Let $C(X)$ denote the algebra of real-valued continuous
functions on $X$. It is clear that every  $\sigma$-$\beta$-defavorable  is 
$\sigma_{C(X)}$-$\beta$-defavorable. We shall show later 
that the converse is no longer true; however, in some  situations
it does as the following statement shows (the straightforward
proof 
is omitted). 
\begin{proposition} Let $X$ be a normal space.  Then  $X$ is
$\sigma_{C(X)}$-$\beta$-defavorable
if and only if $X$ is $\sigma$-$\beta$-defavorable. 
\end{proposition} 
 Christensen-Saint Raymond's game $\mathcal J$ 
was invented  to study the problem of  the existence of continuity points for
separately continuous
mappings. As we shall see,
it is quite possible to replace the game $\mathcal J$
by its variant ${\mathcal J}_\Gamma$ (with suitable
$\Gamma$) and, in this connection,  the next  assertion tells us that 
these games are in a sense the appropriate ones.  
For a set  $\Gamma\subset \mathbb R^X$, let $X_\Gamma$ denote the space
obtained 
when  $X$ is equipped with the topology 
generated by the functions in $\Gamma$ ($C(X_\Gamma)$ stands for the  algebra of
real-valued continuous  functions on  $X_\Gamma$). 
\begin{proposition} Let $X$ be a topological
space and $(K_n)_{n\in\mathbb N}\subset \mathbb R^X$.  Let
$\Gamma=\cup_{n\in\mathbb N}K_n$ and suppose that
 for each $n\in\mathbb N$, the set $A_n$ of $x\in X$ such that
$K_n$ is equicontinuous at $x$ is a residual  subset of $X$. Then  $X$ is
conditionally
$\sigma_{C(X_\Gamma)}^*$-$\alpha$-favorable. In particular, $X$ is
conditionally
$\sigma_{\Gamma}^*$-$\alpha$-favorable.
\end{proposition}
\begin{proof} 
 We shall define a strategy $\tau$ for the Player
 $\alpha$ so that
 for each play which is compatible with $\tau$, say $((V_n,U_n))_{n\in\mathbb N}$,
 the following holds: 
 for every $t\in \cap_{n\in\mathbb N}U_n$, $a_n\in U_n$ (${n\in\mathbb N}$)
 and $g\in \Gamma$, the sequence
 $(g(a_n))_{n\in\mathbb N}$ converges to $g(t)$; in other words,
 the sequence $(a_n)_{n\in\mathbb N}$ converges to $t$ in $X_\Gamma$. Clearly,
such a strategy for $\alpha$ is conditionally winning in the game ${\mathcal
J}_\Gamma^*$. 
 \par
 Let $A=\cap_{n\in\mathbb
N}A_n$ and let us fix  a sequence
$(G_n)_{n\in\mathbb N}$ of dense open subsets of
$X$ such that $\cap_{n\in\mathbb N}G_n\subset A$.
Suppose that $\tau$ has been defined  until stage $n$ and denote by  $V_n$ 
the    $n$th move of Player  $\beta$. Let  ${\mathcal E}_{n}$ be the collection
of all nonempty open sets $U\subset V_n\cap G_n$ 
 such that
$|g(x)-g(y)|\leq 1/n$ for every $x,y\in U$ and $g\in
\cup_{i\leq n}K_i$. Put $\tau(V_n)=V_{n}\cap G_n$ if
${\mathcal E}_{n}$ is empty; if not, choose  $U_{n}\in{\mathcal E}_{n}$
and put  $\tau(V_n)=U_{n}$.
\par  
Let $((V_n,U_n))_{n\in\mathbb N}$ be a play which is compatible with $\tau$, 
 $a_n\in U_n$ ($n\in\mathbb N$),
 $g\in \Gamma$ and  $t\in\cap_{n\in\mathbb N}U_n$.  We have $t\in \cap_{n\in\mathbb
N}G_n$, which implies
 that all the collections ${\mathcal E}_n$, $n\in\mathbb N$, are nonempty. Let
 $p\in\mathbb N$ be  such that $g\in K_p$; since $t,a_n\in U_{n}$, in view of the
choice of the open set $U_n$, we have
 $|g(a_n)-g(t)|<1/n$ for every $n\geq p$. Consequently,
 $\lim g(a_n)=g(t)$.    
 \end{proof}
 The space  $X$ is called  an Eberlein-Grothendieck space (EG-space for short)
if
$X$ is Hausdorff and
 there
 is a compact
 set $\Gamma\subset C(X)$ such that $X=X_\Gamma$. The class of 
  $EG$-spaces includes all metrizable spaces  \cite{arh} (as suggested by the referee, it
  suffices to note that the functions $x\to d(x,y)-d(x_0,y)$, $y\in X$, lie in the
  pointwise compact set of the 1-Lipschitz functions that map the specified point $x_0\in X$ to 0;
  $d$ being a bounded compatible metric on $X$).  Therefore, the next
statement which is a consequence of the proof of Proposition 2.2 improves  the
result
of   Saint
Raymond cited above. 
 \begin{proposition} Let $X$ be an  $EG$-space. Then the following are
equivalent:
 \begin{enumerate}
 \item  $X$ is a Namioka space,
\item $X$ is Baire and conditionally $\sigma$-$\alpha$-favorable,
\item  $X$ is
$\sigma$-$\beta$-defavorable.
\end{enumerate}
 \end{proposition}
 
 \section{Main results}
In what follows, including the statements of Theorems 3.1 and 3.2 and
their respective Corollaries 3.3 and 3.4, 
$f: X\times K\to\mathbb R$ is  a fixed separately continuous mapping and 
  $\phi: K\to
C_p(X)$ is the continuous mapping defined by $\phi(y)(x)=f(x,y)$. We denote
by 
$C_p(X)$  the algebra $C(X)$ equipped with the pointwise
convergence topology.\par
 Let  $\Gamma$ be a nonempty subset
of the product space ${\mathbb R}^X$. A decreasing sequence $(U_n)_{n\in\mathbb
N}$ of  
subsets of  $K$ is said to be {\it countably pair complete with respect to
$(\phi,\Gamma)$} if
for any sequences $(y_n)_{n\in\mathbb N}$, $(z_n)_{n\in\mathbb N}$ such that 
$y_n,z_n\in U_n$
for all
$n\in\mathbb N$,
the sequence $(\phi(y_n)-\phi(z_n))_{n\in\mathbb N}$ has at least one cluster
point
in the subspace  $\Gamma$ of the product space ${\mathbb R}^X$. A sequence  
${\mathcal
U}_n=\{U_k^n: k\in\mathbb N\}$ ($n\in\mathbb N$) of  covers
 of $K$ is said to be  {\it countably pair complete with respect
to $(\phi,\Gamma)$} if for each
 $\sigma\in{\mathbb N}^{\mathbb N}$, the sequence $(\cap_{i\leq
n}U_{\sigma(i)}^i)_{n\in\mathbb
N}$ is countably pair complete with respect to  
 $(\phi,\Gamma)$.
\par
The main results are given in  the next two statements.
The first one should be compared with  \cite[Th\'eor\`eme 5.1]{tal}. The proofs
are postponed  to the next
section.
\begin{theo}  Suppose that there exist a set $\Gamma\subset
\mathbb R^X$ 
 and a sequence   $({\mathcal U}_n)_{n\in\mathbb N}$ of countable covers of
  $K$ such that
 $X$ is $\sigma_\Gamma$-$\beta$-defavorable and    the sequence $({\mathcal
U_n})_{n\in\mathbb N}$ is countably pair complete
 with respect to $(\phi,\Gamma)$.
   Then for every $\varepsilon>0$, there is a residual subset  $R_\varepsilon$
of $X$ such that for every $(x,y)\in R_\varepsilon\times K$ the following holds: 
 \begin{enumerate}
 \item[{$(*)$}] there are  
 a finite sequence $F_i\in{\mathcal U}_i$, $i=0,\ldots, k$, with
$y\in\cap_{i\leq k} F_i$, and  
 a neighborhood $O$ of $(x,y)$ in $X\times K$
 such that: 
 $$|f(x,y)-f(x',y')|<\varepsilon \ {\rm for\,  every}\,  (x',y')\in O\cap
[X\times (\cap_{i\leq k} F_i)].$$ 
\end{enumerate}
\end{theo}
 An important special case of Theorem 3.1 is when 
 $({\mathcal U}_n)_{n\in\mathbb N}$ is a sequence of open covers of the space
$K$; in this case, following a terminology from \cite{kenmoo}, condition $(*)$
says that 
the mapping  $f$ is $\varepsilon$-continuous at the point $(x,y)$ of
$R_\varepsilon\times K$. \par
Recall that a set $A\subset X$ is said to be everywhere of second category 
in $X$ if for every nonempty open set $U\subset X$, the set $A\cap U$ is of the second
category in $U$ (equivalently, in $X$).
 \begin{theo} Let $(U_n)_{n\in\mathbb N}$
 be a sequence of open subsets of $K$, $\Gamma\subset {\mathbb R}^X$ and 
$P$ the set of $y\in K$ for which  there is $\sigma\in{\mathbb N}^{\mathbb N}$
such that $y\in\cap_{n\in\mathbb N}U_{\sigma(n)}$ and the sequence
$(\cap_{i\leq n}U_{\sigma(i)})_{n\in\mathbb N}$ is countably pair complete
with respect to  $(\phi,\Gamma)$. Denote by $R_\varepsilon(P)$
the set of  $x\in X$
such that  the mapping $f: X\times K\to\mathbb R$ is $\varepsilon$-continuous
at each point of $\{x\}\times P$.  Then
\begin{enumerate} 
\item if $X$ is conditionally $\sigma_\Gamma$-$\alpha$-favorable, 
$R_\varepsilon(P)$ is
a residual subset of $X$;
\item if $X$ is  $\sigma_\Gamma$-$\beta$-defavorable,
$R_\varepsilon(P)$ is everywhere of second
category in $X$.
\end{enumerate}
\end{theo}
 To express some consequences
 of these results, we need to recall some terminology. Let $Z$
 be a topological space and $Y\subset Z$. The set $Y$ is said to be bounded (or
 relatively pseudocompact) in $Z$
if every continuous function $g:Z\to\mathbb R$ is bounded on $Y$;
 $Z$ is pseudocompact if $Z$ is Tychonoff (i.e., completely regular)  and
bounded in itself . The space
$Y$ is called  Lindel\"of in $Z$ if 
every open cover of $Z$ has a countable subcover of $Y$.
The space $Y$ is said to be 
relatively  countably compact in $Z$ if  every sequence $(y_n)_{n\in\mathbb
N}\subset Y$
has a cluster point  in $Z$. \par

Also let us remind  that a Tychonoff space $Y$ is called a  $p$-space if
there is 
a sequence  $({\mathcal U}_n)_{n\in\mathbb N}$ of open covers of $Y$ such that
 to each   $x\in Y$ corresponds a sequence $U_n\in{\mathcal U}_n$, $n\in\mathbb
N$,
 such that  $x\in\cap_{n\in \mathbb N}U_n$ and the intersection
$\cap_{n\in\mathbb N}U_n$
 is a compact subset of $Y$ for which the sequence $(\cap_{i\leq
n}U_i)_{n\in\mathbb N}$
 is an outer basis. Finally,
a Tychonoff space $Y$ is said to be  \v Cech-complete if there
is a sequence   $({\mathcal U}_n)_{n\in\mathbb N}$ of open covers of  $Y$ which
is complete
in the sense that  any closed filter basis  ${\mathcal F}$ on $Y$
has a nonempty intersection, provided that for each 
 $n\in{\mathbb N}$ there are $F\in{\mathcal F}$ and
$U\in{\mathcal U}_n$ such that $F\subset  U$.
 
\begin{cor} Suppose that $X$ is a $\sigma_{C(X)}$-$\beta$-defavorable space.
Then, there is a $G_\delta$ dense subset $R$ of $X$ such that $f$ is jointly
continuous at each point of $R\times K$, in each of the following cases:
\begin{enumerate}
\item $K$ is  pseudocompact  and  every bounded subspace of $C_p(X)$
  is relatively countably compact in $C_p(X)$.
\item $K\times K$ is pseudocompact and  every pseudocompact subspace of $C_p(X)$
  is relatively countably compact in $C_p(X)$.
\item  $\phi(K)$ is relatively  Lindel\"of in a \v Cech-complete
subspace of
$ C_p(X)$.
\item  $K$ is  Lindel\"of \v Cech-complete.
\end{enumerate}
\end{cor}
\begin{proof} We apply Theorem 3.1 by taking   $\Gamma=C(X)$ in each of these
cases. For (1) and (2), let ${\mathcal U}_n=\{K\}$ for every
$n\in\mathbb N$,  and note that  the
set $L=\phi(K)-\phi(K)$ is bounded
in $C_p(X)$. Indeed, $L$ is  a pseudocompact
subspace of $C_p(X)$ in case (2); in case (1), the set $L$ is the difference
of two pseudocompact subspaces of the topological group
$C_p(X)$, hence, 
according to a result of Tka\v cenko \cite{tka}, it is
bounded in $C_p(X)$. \par 
For (3), let $Z$ be a \v Cech-complete subspace of $C_p(X)$ such that $\phi(K)$
is Lindel\"of in $Z$.  Let $({\mathcal W}_n)_{n\in\mathbb N}$ be a complete
sequence of open covers of $Z$; since $\phi(K)$ is Lindel\"of in $Z$, for each
$n\in\mathbb N$ there is a countable collection ${\mathcal V}_n\subset{\mathcal
W}_n$ such that $\phi(K)\subset\cup{\mathcal V}_n$. The
sequence $(\phi^{-1}({\mathcal V}_n))_{n\in\mathbb N}$ fulfills 
the conditions of  Theorem 3.1.\par
The proof of (4) is similar to (3).
\end{proof}

The point (4) in Corollary 3.3 is established in \cite{kenmoo} for $X$
Lindel\"of
$\alpha$-favorable. The point (1) in the following is proved 
in \cite{moo} for $X$
conditionally $\sigma$-$\alpha$-favorable;  the point
(2) describes the situation in the ``$\beta$-defavorable'' case. 
\begin{cor} Suppose that  $K$ is a Lindel\"of $p$-space and let $\varepsilon>0$.
Let $R_\varepsilon$
be the set of $x\in X$ such that $f$ is $\varepsilon$-continuous
at each point of $\{x\}\times K$.
\begin{enumerate} 
\item If $X$ is  conditionally $\sigma_{C(X)}$-$\alpha$-favorable,
then $R_\varepsilon$ is a residual subset of $X$.
\item If $X$ is $\sigma_{C(X)}$-$\beta$-defavorable, then $R_\varepsilon$
is everywhere of second category in $X$.
\end{enumerate}
\end{cor}
\begin{proof}  Since $K$ is a Lindel\"of $p$-space, letting  $P=K$, Theorem 3.2
applies.
\end{proof}

\section{The proofs}
Lemma 4.1 below is used repeatedly in the  proof   
of Theorem 3.1. The proof of  Theorem 3.2 consists in  adapting
that of Theorem 3.1;
 Lemma 4.1 is not needed there, however, for the first item, it  is replaced
by the well-known characterization
of residual sets in term of the Banach-Mazur game (see below).
So we give the entire proof for Theorem 3.1 and only indicate the main changes
to get Theorem 3.2. \par
Lemma 4.1 is an immediate
 consequence of the following well-known property  \cite{kel}: ``given a set 
$A\subset X$ 
which is of  second category in    $X$, there is 
 a nonempty open subspace $V$ of $X$ such that   $A\cap V$ is everywhere  of 
second category in   $V$''. (The concept is recalled just before Theorem 3.2).  

 \begin{lemma} Let $A=\cup_{n\in\mathbb N}B_n$ be a set  of
the second
 category in the space
 $Y$. Then, there is a nonempty open set $V\subset Y$ and $n\in\mathbb N$ such
that
 $B_n\cap V$ is everywhere of  second category in $V$.
 \end{lemma}

\vskip 2mm
\noindent {\bf Proof of   Theorem 3.1} As said in the introduction, the main
arguments  follow the proof given--in French-- in \cite{bo1} for 
Christensen-Saint Raymond's theorem cited above.\par  
\vskip 2mm
\noindent {\it The assumption.}   For  $\varepsilon>0$, $F\subset X$ and
$L\subset K$, let
$R_\varepsilon(F,L)$ (or simply
  $R(F,L)$) be the set of  $x\in F$ such that the property $(*)$ is satisfied
for all $y\in L$. We have 
  to prove that
 $R(X,K)$ is a residual subset of $X$. Let us suppose to the contrary and  show
 that $X$ is $\sigma_{\Gamma}$-$\beta$-favorable.  Thus, writing 
$D(F,L)=F\setminus R(F,L)$, our assumption says that the set $D(X,K)$
is of  second category in $X$.\par 
 \vskip 2mm 
  \noindent {\it The  strategy}.  Write 
  ${\mathcal U}_n=\{F_k^n:k\in\mathbb N\}$. 
 We are going to define a strategy $\sigma$
for the Player $\beta$ in the game ${\mathcal J}_\Gamma$  which produces parallel  
to
each
play $(V_n, (a_n,U_n))_{n\in\mathbb N}$ 
a set of sequences $(x_n)_{n\in\mathbb N}, (t_n)_{n\in\mathbb N}\subset X$,
$(y_n)_{n\in\mathbb N}, (z_n)_{n\in\mathbb N}\subset K$ and $(k_n)_{n\in\mathbb
N}\subset \mathbb N$, 
so that for every $n\in\mathbb N$:
\begin{enumerate}
\item the set $D(V_n,\cap_{i\leq n}F_{k_i}^i)$ is everywhere 
of the second category in  $V_n$;
\item  $y_{n+1},z_{n+1}\in\cap_{i\leq n}F_{k_i}^i$;
\item $|f(a_i,z_{n+1})-f(a_i,y_{n+1})|< 1/(n+1)$ for each $0\leq i\leq n$;
\item $V_{n+1}\subset\{t\in X:
|f(t,z_{n+1})-f(t_{n+1},z_{n+1})|<\varepsilon/3\}\cap
\{t\in X: |f(t,y_{n+1})-f(x_{n+1},y_{n+1})|<\varepsilon/3\}$;
\item $|f(x_{n+1},y_{n+1})-f(t_{n+1},z_{n+1})|\geq \varepsilon$.

\end{enumerate}

Applying
 Lemma 4.1  to  $Y=X$ and $A=D(X,K)$  gives
 a nonempty open set  $V_0\subset X$  and $k_0\in\mathbb N$ such that
 $D(V_0,F_{k_0}^0)$ is everywhere of second category
in $V_0$.  Let $y_0,z_0\in F_{k_0}^0$, $x_0,t_0\in
X$
be arbitrary  and define $\sigma(\emptyset)=V_0$. Assume that we are at stage 
$p$:
Player  $\alpha$ having produced
$(a_0,U_0),\ldots, (a_p, U_p)$, Player  $\beta$ his sequence
$V_0,\ldots,V_p$ and all  terms of sequences above having been defined until $p$
in accordance   
with  (1)-(5). First, let $x_{p+1}\in U_p$
and 
$y_{p+1}\in \cap_{i\leq p}F_{k_i}^i$ be so that the condition $(*)$ is not
satisfied (the inductive hypothesis (1)
ensures that  $D(U_p, \cap_{i\leq p}F_{k_i}^i)$
is not empty). The set $$A=\{t\in U_p:
|f(t,y_{p+1})-f(x_{p+1},y_{p+1})|<\varepsilon/3\}$$
is a neighborhood  of $x_{p+1}$ in $X$ and the set
$$B=\cap_{i\leq p}\{z\in K:
  |f(a_i,z)-f(a_i,y_{p+1})|<\varepsilon/4\}$$ 
is a neighborhood of  $y_{p+1}$ in $K$; choose
$(t_{p+1},z_{p+1})\in A\times [B\cap(\cap_{i\leq p}F_{k_i}^i)]$ such that
$|f(t_{p+1},z_{p+1})-f(x_{p+1},y_{p+1})|\geq\varepsilon$. The open set
$$ O=A\cap\{t\in U_p: |f(t,z_{p+1})-f(t_{p+1},z_{p+1})|<\varepsilon/3\}$$
being nonempty $(t_{p+1}\in O$), the set  $D(O,\cap_{i\leq p}F_{k_i}^i)$ is of
the second category
in
$O$; since  $F_{k_p}^p\subset\cup_{l\in \mathbb N}F_l^{p+1}$,  Lemma 4.1
gives
an integer $k_{p+1}$ and a nonempty open set $V_{p+1}\subset O$ such that
  $D(V_{p+1},\cap_{i\leq p+1}F_{k_i}^i)$ is everywhere of second category in 
  $V_{p+1}$. Define
$$\tau((a_1,U_1),\ldots, (a_p, U_p))=V_{p+1}.$$
All items (1)-(5) are satisfied  for $i\leq p+1$.  The definition
of the strategy  $\sigma$ is complete.\par
\vskip 2mm
\noindent {\it Conclusion.} We show that $\sigma$
is a winning strategy for Player $\beta$.
Suppose that    $((a_n,U_n))_{n\in\mathbb N}$  is a winning  play for
 $\alpha$ against the strategy $\sigma$. According to (2), there is a cluster
 point $g\in\Gamma$ 
 of the
sequence $(\phi(z_n)-\phi(y_n))_{n\in\mathbb N}$. 
 According to (3), for every
$m\in\mathbb N$, we
 have $$\lim_n|f(a_m,z_n)-f(a_m,y_n)|=0;$$ thus $g(a_m)=0$ for every
$m\in\mathbb N$. It follows that there is
  $t\in
\cap_{n\in\mathbb
N}U_n$  such that $g(t)=0$; in particular,
$|f(t,z_n)-f(t,y_n)|<\varepsilon/3 $ for some $n\in\mathbb N$. It follows from
 (4) that
\begin{align}
\nonumber
|f(x_n,y_n)-f(t_n,z_n)|&\leq  |f(x_n,y_n)-f(t,y_n)|+|f(t,y_n)-f(t,z_n)|\\
\nonumber
&\ \ +  |f(t,z_n)-f(t_n,z_n)|\\ \nonumber
&<  \varepsilon,
\end{align}
contrary to  (5). 
\begin{remark} {\rm Suppose that $X$ is  $\sigma$-$\beta$-defavorable. Then,
the argument in the conclusion
step of  the above proof also works if the assumption on $K$ is weakened 
assuming that
the sequence   $({\mathcal U}_n)_{n\in \mathbb N}$ is countably pseudo-complete
with respect
to $\phi$ in
the  following sense:
For every  $F_n\in{\mathcal U}_n$, $n\in\mathbb N$, and every sequence
$(y_n)_{n\in\mathbb N}\subset K$
such that $y_n\in\cap_{i\leq n}F_i$, the set $\{\phi(y_n):n\in\mathbb N\}$
is bounded in $C_p(X)$.  Indeed, let $((a_n,U_n))_{n\in\mathbb N}$  be a winning
 play for
 $\alpha$ against the strategy $\sigma$ and choose $t\in\cap_{n\in\mathbb N}U_n$
 such that $t\in\overline{\{a_n:n\in\mathbb N\}}$. Write
$A=\{t\}\cup\{a_n:n\in\mathbb N\}$
 and let  $r_A: C_p(X)\to C_p(A)$ be the map under which  each $g\in C(X)$ is
sent to its restriction to $A$.
The sets $\{r_{A}(\phi(y_n)):n\in\mathbb N\}$ and
$\{r_{A}(\phi(z_n)):n\in\mathbb N\}$
are bounded thus relatively compact in $C_p(A)$, since $C_p(A)$ is metrizable
(see for instance
Lemma III.4.7 in \cite{arh}); it follows that  the sequence
$\big(r_A(\phi(z_n))-r_A(\phi(y_n))\big)_{n\in\mathbb N}$ has a cluster  point $g\in
C_p(A)$. By (3),
$g(\{a_n:n\in\mathbb N\})=\{0\}$, hence $g(t)=0$ and  the proof can be continued
as above. 
Let us mention that the corresponding result (that is, the property ${\mathcal
N}(X,K)$ holds) in case when $\phi(K)$ is bounded in $C_p(X)$ (and $X$ is
$\sigma$-$\beta$-defavorable) is due to
Troallic  \cite{tro}. Unfortunately, there is no hope to establish the
same result in case  when $X$
is $\sigma_{C(X)}$-$\beta$-defavorable (see Example 5.3 below).} 
\end{remark}
\vskip 2mm
Before we pass to  Theorem 3.2, let us recall the description
of first category sets in term of the Banach-Mazur game. For
 a space $Y$ and  $R\subset  Y$, a play in the game  $BM(R)$ (on $Y$) 
is a sequence  $((V_n,U_n))_{n\in\mathbb N}$ of pairs of nonempty open subsets
of $Y$ produced  alternately by two players $\beta$
and $\alpha$ as  follows:  $\beta$ is the first to move and gives $V_0$, then
Player $\alpha$ gives $U_0\subset V_0$; at stage $n\geq 1$,
the open set $V_n\subset U_n$ being chosen by $\beta$, Player $\alpha$ gives
$U_n\subset
V_n$. 
Player $\alpha$ wins the play if
 $\cap_{n\in\mathbb N}U_n\subset R$. It is well known that
$X$ is $BM(R)$-$\alpha$-favorable (i.e., $\alpha$
has a winning strategy in the game $BM(R)$)  if and only if  $R$
is a residual subset of $Y$. The reader is referred to  \cite{oxt}.
\vskip 2mm
\noindent{\bf Proof of  Theorem 3.2.} Denote by ${\mathbb N}^{<{\mathbb N}}$ the
set
of finite sequences of integers and
let $\phi: {\mathbb N}^{<{\mathbb N}}\to\mathbb N$ be a bijective map such that
$\phi(s)\geq |s|$
for every $s\in {\mathbb N}^{<{\mathbb N}}$, where $|s|$
stands for the length of $s$.  For 
$n\in\mathbb N$, define 
$$F_n=\bigcap_{i\leq|\phi^{-1}(n)|}U_{\phi^{-1}(n)(i)}.$$
We keep the notation $D(F,L)$ (for $F\subset X$ and $L\subset K$) used
in the proof of Theorem 3.1 and write $R(P)$ for $R_\varepsilon(P)$.\par 
1) Let $\tau_1$ be a conditionally winning
strategy for Player $\alpha$
in the game ${\mathcal J}_\Gamma$. We deduce from $\tau_1$ a
winning strategy $\tau_2$ for Player  $\alpha$ in the game $BM(R(P))$ as
follows.
Fix $*\not\in K\cup X$. Let  $V_n$ be the $n$th move   of $\beta$ 
in the game $BM(R(P))$ and write $(W_n, a_n)=\tau_1(V_0,\ldots, V_n)$. If
$D(W_n,F_n)=\emptyset$, define $\tau_2(V_0,\ldots, V_n)=W_n$ and 
 $x_n=t_n=y_n=z_n=*$. If
 $D(W_n,F_n)\not=\emptyset$,  first choose $x_n\in W_n$ and $y_n\in F_n$
such that $f$ is not $\varepsilon$-continuous at the point $(x_n,y_n)$.
Then, considering the sets
$$A=\{t\in W_n: |f(t,y_n)-f(x_n, y_n)|<\varepsilon/3\}$$ and
$$B=\cap_{i\leq n}\{z\in K: |f(a_i,z)-f(a_i, y_n)|<1/(n+1)\},$$
choose $t_n\in A$ and $z_n\in B\cap F_n$ such that 
$|f(x_n,y_n)-f(t_n, z_n)|\geq\varepsilon$; finally define

$$\tau_2 (V_0,\ldots, V_n)=\{t\in W_n:|f(t,z_n)-f(t_n, z_n)|<\varepsilon/3\}.$$
The definition of $\tau_2$ is complete.\par 
Let us suppose for contradiction that  Player $\beta$
 has a winning play 
$(V_n)_{n\in\mathbb N}$ against 
the strategy $\tau_2$. Then 
$\cap_{n\in\mathbb
N}V_n\not\subset R(P)$, that is, there are $a\in\cap_{n\in\mathbb N}V_n$ and
$y\in P$
such that $f$ is not $\varepsilon$-continuous at $(a,y)$. Let
$\sigma\in{\mathbb N}^{\mathbb N}$ be such that the
sequence $(\cap_{i\leq n}U_{\sigma(i)})_{n\in\mathbb N}$
is countably pair complete with respect
to $(\phi,\Gamma)$ and $y\in\cap_{i\in\mathbb N} U_{\sigma(i)}$.
For  $n\in\mathbb N$, 
let $k_n=\phi(\sigma_{|n})$; then  $a\in W_{k_n}$ and $y\in F_{k_n}$,
thus $D(W_{k_n},F_{k_n})\not=\emptyset$ which indicates that  $y_{k_n},z_{k_n}$
have been selected
in $F_{k_n}$.  Since  $F_{k_n}=\cap_{i\leq n}U_{\sigma(i)}$, the sequence
$(\phi(y_{k_n})-\phi(z_{k_n}))_{n\in\mathbb N}$ has at least
a cluster point   $g\in \Gamma$. Since $\cap_{n\in\mathbb N}V_n\not=\emptyset$
and
the  strategy $\tau_1$ is  conditionally winning, there is  $t\in
\cap_{n\in\mathbb N}V_n$ such that $g(t)\in\overline{\{g(a_n):n\in\mathbb N\}}$
(note that the play $(V_n,(W_n,a_n))_{n\in\mathbb N}$ is compatible with
$\tau_1$). The  argument    
from the ``conclusion" step in the proof of Theorem 3.1 can now be used
to get the required contradiction. Therefore,
$R(P)$  is a residual subset of $X$.\par
2) We proceed as in (1), keeping the same notations.  Suppose
that there exists a nonempty open set $\Omega$ such that $R(P)\cap
\Omega$
is of
first category in $X$, that is, $R(P)\cap
\Omega\subset\cup_{n\in\mathbb 
N}A_n$ where each $A_n$ is a closed nowhere dense subset of $X$.
We deduce from this
a winning  strategy $\sigma$ for Player $\beta$ 
in the game ${\mathcal J}_\Gamma$ as follows.  To begin let
$\sigma(\emptyset)=\Omega$. At step $n$,
let   $(V_0,a_0),\ldots, (V_n,a_n)$ be the first $n$th moves of Player $\alpha$
and
consider the nonempty open set $O_n=V_n\setminus A_n$. If
$D(O_n,F_n)=\emptyset$,
define $\sigma((V_0,a_0),\ldots, (V_n,a_n))=O_n$ and
$t_n=x_n=y_n=z_n=*$;  if  $D(O_n,F_n)\not=\emptyset$,
define

$$\sigma((V_0,a_0),\ldots, (V_n,a_n))=\{t\in O_n: |d(f(t,z_n),f(t_n,
z_n)|<\varepsilon\},$$
 the  points $x_n$, $t_n$, $y_n$, $z_n$ being chosen exactly  as in (1).\par 
Suppose that
$((a_n,V_n))_{n\in\mathbb
N}$
is a  play for Player $\alpha$ which is compatible with $\sigma$ and let us show
that there is $g\in \Gamma$ so that $g(t)\not\in\overline{\{g(a_n):n\in\mathbb
N\}}$
for every
 $t\in \cap_{n\in\mathbb N}V_n$. Since $\Gamma\not=\emptyset$,  we may assume
that $\cap_{n\in\mathbb
N}V_n\not=\emptyset$. 
Let $a\in\cap_{n\in\mathbb N}V_n$; then $a\not\in 
R(P)$ hence there is
$y\in P$ such that
 $f$ is not $\varepsilon$-continuous at the point $(a,y)$.
Since
$y\in P$, there is 
$\sigma\in{\mathbb N}^{\mathbb N}$ such that $y\in\cap_{n\in\mathbb
N}U_{\sigma(n)}$ and the
sequence
$(\cap_{i\leq n}U_{\sigma(i)})_{n\in\mathbb N}$ is countably pair 
complete with respect to $(\phi,\Gamma)$.  
As in (1), letting 
 $k_n=\phi(\sigma_{|n})$ for  $n\in\mathbb N$, we obtain that   $a\in V_{k_n}$
and $y\in
F_{k_n}$,
hence $D(O_{k_n},F_{k_n})\not=\emptyset$ and,   consequently,
$\{y_{k_n},z_{k_n}\}\subset F_{k_n}$ for every
$n\in\mathbb N$. Take  a cluster point $g\in \Gamma$ of the sequence 
$(\phi(y_{k_n})-\phi(z_{k_n}))_{n\in\mathbb N}$; the assumption
that $g(t)\in\overline{\{g(a_n):n\in\mathbb
N\}}$ for some
$t\in\cap_{n\in\mathbb
N}V_n$ leads  to a contradiction as in (1).  
\vskip 2mm
To conclude this section, let us mention the following result which
gives  a  description of the class  of Namioka spaces and answers in a certain sense Question 1167 (or Question 8.2) in \cite{bormoo}.
\begin{proposition} Let  $X$ be a Tychonoff  space. Then the following are
equivalent:
\begin{enumerate}
\item $X$ is a Namioka space,
\item $X$ is a Baire space and conditionally
$\sigma_{\Gamma}^*$-$\alpha$-favorable
for every compact $\Gamma\subset C_p(X)$,
\item $X$ is a Baire space and  conditionally
$\sigma_\Gamma$-$\alpha$-favorable for every compact $\Gamma\subset C_p(X)$,
\item   $X$  is $\sigma_\Gamma$-$\beta$-defavorable
for every compact $\Gamma\subset C_p(X)$.
\end{enumerate}
 \end{proposition}
 \begin{proof} The fact that (1) implies (2) follows from Proposition 2.2 and Saint Raymùond's theorem
 that every Tychonoff Namioka space is Baire \cite{sr2}. The implications
 $(2)\to (3)$ and $(3)\to (4)$ are obvious. Finally, Theorem 3.1 shows that (4) implies (1). 
 \end{proof}
\section{Some  related results}
Recall that a subspace   $X$ of a topological space $Y$ is said to be
$C$-embedded
in $Y$ if 
every  $f\in C(X)$ has an extension $g\in C(Y)$. Suppose that $X$
is dense
in $Y$; it is well known that  $X$ is  $C$-embedded in $Y$ if and only if
$X$ is  $G_\delta$-dense in  $Y$ and $z$-embedded in  $X$, that is,  every
zero set
of $X$ is the intersection with $X$ of a zero set of $Y$. \par 

To establish the following proposition we note that the rule that  Player
$\alpha$ wins the play $((U_n, V_n))_{n\in\mathbb N}$ 
in the game ${\mathcal J}_{C(X)}^*$ (the strong version of ${\mathcal 
J}_{C(X)}$) can be formulated in an equivalent manner as follows: For every 
zero set $Z\subset X$ such 
that $Z\cap U_n\not=\emptyset$
for every $n\in\mathbb N$, we have $Z\cap(\cap_{n\in\mathbb
N}U_n)\not=\emptyset$.
 
\begin{proposition} Let $Y$ be a  $\sigma_{C(Y)}^*$-$\beta$-defavorable space
(respectively, $\sigma_{C(Y)}^*$-$\alpha$-favorable).
 Then every
$C$-embedded dense subspace $X$ of $Y$   
is $\sigma_{C(X)}^*$-$\beta$-defavorable (respectively,
$\sigma_{C(X)}^*$-$\alpha$-favorable).
\end{proposition}
\begin{proof} We outline a proof of  the $\beta$-defavorable case (the other case
being
similar). Let $\tau_X$
be a strategy for Player  $\beta$ in the game ${\mathcal J}_{C(X)}^*$ and let us
show that it
is not a winning one. Fix a map $V\to V^*$ under which   each nonempty open subset of $X$ 
is sent to an open
subset $V^*$ of $Y$ such that $V=V^*\cap X$.
Consider the following strategy  $\tau_Y$ for Player $\beta$ in the game
${\mathcal J}_{C(Y)}^*$. Write $V_0=\tau_X(\emptyset)$ and 
put $\tau_Y(\emptyset)=V_0^*$.   Suppose that $\tau_Y$
has been defined until stage $n$ and write $V_{n+1}=\tau_X(U_0\cap X,\ldots,U_n\cap X)$,
where $U_0,\ldots, U_n$ are the first $n+1$ moves of Player  $\alpha$ in  the game
${\mathcal J}_{C(Y)}^*$. Define $\tau_Y(U_0,\ldots, U_n)=V_{n+1}^*\cap U_n$ 
(this open subset of $U_n$ is nonempty 
because it contains $V_{n+1}$).\par 
 There is a winning play $(U_n)_{n\in\mathbb N}$
for Player $\alpha$ against the strategy $\tau_Y$ in the game
${\mathcal J}_{C(Y)}^*$. The 
corresponding sequence $(U_n\cap X)_{n\in\mathbb N}$ is a play with respect to the
game  ${\mathcal J}_{C(X)}^*$, which is compatible   with
 $\tau_X$.
Let 
$Z$ be a  zero set of  $X$ such that $Z\cap U_n\not=\emptyset$
for every $n\in\mathbb N$. There is  a zero set   $T$  of $Y$ such that $Z=T\cap
X$;
the set $H=T\cap(\cap_{n\in\mathbb N}U_n)$ is a nonempty  $G_\delta$
subset of $Y$;  since   $X$ is $G_\delta$-dense
in  $Y$, we obtain $Z\cap(\cap_{n\in\mathbb N}U_n)=H\cap X\not=\emptyset$. 
\end{proof}

A standard example illustrating Proposition 5.1 is when
  $X$ is pseudocompact   and $Y$ is its  Stone-\v Cech-compactification $\beta
X$. Clearly,  $\beta X$ (as any compact space)
is $\sigma_{C(\beta X)}^*$-$\alpha$-favorable; thus Proposition 5.1 leads to
the following.

\begin{cor} Every   pseudocompact space $X$ is
$\sigma_{C(X)}^*$-$\alpha$-favorable.
\end{cor}
It follows from Corollary 5.2 and Proposition 4.3 that
every pseudocompact space is a Namioka space. Actually
a stronger statement can be established (see Proposition 5.5 below).
We are now ready to give an example of a 
$\sigma_{C(X)}^*$-$\alpha$-favorable, hence
$\sigma_{C(X)}$-$\beta$-defavorable, which is not
$\sigma$-$\beta$-defavorable.
\begin{example} {\rm It is shown by Shakhmatov in  \cite{sha} that  there exists
a pseudocompact space $P$ without isolated points, every
countable subset of which
  is discrete. Such a space is  $\sigma_{C(P)}^*$-$\alpha$-favorable in
view of
Corollary 5.2. Using the fact that
$P$ has no isolated point and all its countable subspaces are closed, it is easy
to check
that $P$
is $\sigma$-$\beta$-favorable.
We have mentioned in Remark 4.2 that the property ${\mathcal N}(X,K)$ is
generally false  if $K$ is  pseudocompact and $X$ is
 $\sigma_{C(X)}$-$\beta$-defavorable. Indeed, Shakhmatov's
space $P$ is such that the unit ball $K$ of $C_p(P)$  is
pseudocompact (see for instance Example I.2.5 in \cite{arh} or \cite{tk1}) and
since $P$ has no isolated
point, the evaluation mapping $e:P\times K\to [0,1]$ does not have any point of
continuity.}
 \end{example}
\begin{proposition} Let $X$ be a space such that every countable subspace
of $X$ is  $C$-embedded in $X$. Then,
the space $Y=C_p(X)$ is  $\sigma_{C(Y)}^*$-$\alpha$-favorable.
\end{proposition}
\begin{proof} 
For each cardinal number $\gamma$, the product space $\mathbb R^\gamma$ is 
$\sigma_{C({\mathbb
R}^\gamma)}^*$-$\alpha$-favorable; we refer to
Christensen's paper \cite{chr} for a similar result about the product of
$\tau$-well $\alpha$-favorable spaces (defined therein).
Let
$\nu Y$ stand for  the realcompactification
of $Y$;  then $\nu Y=\mathbb R^X$ \cite{tk1}. Since $Y$ is
$C$-embedded  in  $\nu Y$,
Proposition 5.1 shows that $Y$ is 
$\sigma_{C(Y)}^*$-$\alpha$-favorable.
\end{proof}  
We should  conclude  under the assumption of Proposition 5.4 that the space
$C_p(X)$ is a Namioka space, but Corollary 5.7 below provides
a more general statement.
\begin{proposition} Let $X$ be a Baire space with a dense 
 $\sigma$-bounded subspace. Then, $X$  is a  Namioka space. 
\end{proposition}
\begin{proof} 
Let  $\Gamma$ be a compact subset of  $C_p(X)$ and let us show that
$X$ is $\sigma_{\Gamma}$-$\beta$-defavorable. Recall that
every compact space 
$L$ such that  $C_p(L)$ contains  a $\sigma$-compact subset separating the
points
of $L$ is an Eberlein compactum (\cite{arh}, p. 124). Let  $e: X\to C_p(\Gamma)$
be the mapping
 $e(x)(y)=y(x)$. Since $e$ is  continuous and 
every bounded subset of $C_p(\Gamma)$ is  relatively compact (by the
generalization of Grothendieck's
theorem in \cite{arh}), the closure of $e(X)$  in $C_p(\Gamma)$
contains a 
dense  $\sigma$-compact space $Y$. Clearly, $Y$
separates the points of $\Gamma$,  hence $\Gamma$ is an Eberlein compactum.
By a result of Deville \cite{dev}, every Eberlein compactum is  co-Namioka;
thus,
following
Proposition 2.2, the space $X$
is  $\sigma_\Gamma$-$\beta$-defavorable.
\end{proof}
\begin{remark} {\rm Following \cite{arh}, a space $X$ is called  $k$-primary
Lindel\"of if
$X$ is the continuous image of a closed subspace
of a space of the form $K\times (L(\gamma))^\omega$, where $K$
is a compact space and  $\gamma$ is cardinal number; $L(\gamma)$ stands for the
one point
Lindel\"ofication of the discrete space of cardinality $\gamma$. As suspected in
\cite{merneg}, Remark 2.17, it can be proved that
every  Baire space with a dense
$k$-primary Lindel\"of subspace is a Namioka space. This can
be established with  the same method as in the proof of Proposition 5.5,
replacing Deville's result by  Debs's theorem  that every Corson compactum
is co-Namioka \cite{de2}, and using 
the following theorem by Bandlow \cite{bad}: If $X$ is a $k$-primary Lindel\"of
space,
then every compact subspace of $C_p(X)$ is a Corson compactum.}
\end{remark}
Recall that a space $X$ is called $b$-discrete if every countable subspace $A$
of $X$ is discrete and $C^*$-embedded in $X$ (every bounded continuous function
on $A$ has a continuous extension over $X$). 
\begin{cor} Let  $X$ be a space such that $C_p(X)$ is Baire. If
$X$ is   $b$-discrete, then $C_p(X)$ is a Namioka space.
\end{cor}
\begin{proof} Since $X$ is $b$-discrete, the subspace  
$C^*(X)\subset C_p(X)$ of bounded continuous functions
is  $\sigma$-bounded \cite{tk1}. Since  $C^*(X)$ is dense in $C_p(X)$,
Proposition 5.5 applies.
\end{proof}
The converse of 5.7 is not true as the following example shows.  
\begin{example} {\rm Example 7.2 in \cite{lutmcc} exhibits  a countable   space
$X$ containing
a non-$C^*$-embedded subspace, such that
$C_p(X)$ is Baire. Since
 $C_p(X)$ is metrizable (and Baire) it is  a
Namioka space by the result of Saint Raymond mentioned in the introduction.}
\end{example}
In view of Corollary 5.7 and Example 5.8, it seems likely that
the space $C_p(X)$ (for a Tychonoff space $X$) is a Namioka space
as soon as it is Baire.

\begin{example}{\rm  There is a  Namioka space $X$ which is 
$\sigma_{C(X)}$-$\beta$-favorable. (This is related to Proposition 4.3.)
We give two examples of such  spaces. \par 
1) Let $X$ be the reals equipped with the so-called density topology $T_d$
\cite{oxt}. The space
$X$
is a Namioka space, because it is a Baire space \cite{lmz} and every compact
subset of $C_p(X)$
is metrizable (see \cite{fre} for a general statement). To show that
$X$ is $\sigma_{C(X)}$-$\beta$-favorable, consider  the strategy $\tau$ for
Player
$\beta$ defined as follows: $\tau(\emptyset)=X$ and $\tau((a_0,U_0),\ldots,
(a_n, U_n))=V_{n+1}$, where $V_{n+1}$ is a nonempty open subset of $U_n$
such that $V_{n+1}\subset U_n\setminus\{a_0,\ldots,a_n\}$ and $|x-y|\leq 1/n$
for each $x,y\in V_n$ (recall that $T_d$ is finer than the usual topology).
Suppose that $((V_n,U_n,a_n))_{n\in\mathbb N}$ is
a play which is compatible with $\tau$. It is well known
that every countable subset of $X$ is closed \cite{oxt}; thus, since the
intersection
 $A=\cap_{n\in\mathbb N}V_n$ contains at most one point  (and $X$ is Tychonoff),
there is a function $f\in C(X)$ such that $f_{|A}=0$ and $f(a_n)=1$
for every $n\in \mathbb N$. Thus $\tau$ is a winning strategy. \par
2) If the Continuum Hypothesis (CH) is assumed then there is a Namioka
space $X$ and a countably compact subspace
$\Gamma$ of $C_p(X)$ such that $X$ is $\sigma_\Gamma$-$\beta$-favorable.
Namely, under CH, Burke
and  Pol proved in   \cite{burpol} that the product  $B=\{0,1\}^{\aleph_1}$
equipped
with the so-called Baire topology, that is, the  $G_\delta$-modification
of the usual product $\{0,1\}^{\aleph_1}$, is  a Namioka space. The subspace
 $\Gamma=\{f\in C(B):f(B)\subset\{0,1\}\}$ of  
 $C_p(B)$ is $\omega$-compact, i.e., every countable subset of $\Gamma$
 is relatively compact in $\Gamma$. (See \cite{burpol} or use
Arhangel'ski\v \i 's result that for every  $P$-space $Y$, the space
$C_p(Y,[0,1])$ is  $\omega$-compact \cite{arh}.) A winning strategy $\tau$ for
Player
 $\beta$ in the game ${\mathcal  J}_\Gamma$  consists of producing clopen sets
such that  
 $\tau((a_0,U_0),\ldots, (a_n, U_n))\cap\{a_i:i\leq n\}=\emptyset$ (where
$(U_0,a_0),\ldots, (U_n,a_n)$
 are the first $n$th moves of Player $\alpha$). 
Such a  strategy 
 is indeed  winning  for 
 if $(V_n,(U_n,a_n))_{n\in\mathbb N}$ is a compatible play, then 
 the sequence $(1_{V_n})_{n\in\mathbb N}\subset \Gamma$ has a cluster point 
 $f\in \Gamma$ (in fact, $(1_{V_n})_{n\in\mathbb N}$ converges to
 $1_{\cap_{n\in\mathbb N}V_n}$). Then, since  $f(a_n)=0$ for each $n\in\mathbb
N$, there is no
 point $t\in\cap_{n\in\mathbb N}V_n$ for which 
$f(t)\in\overline{\{f(a_n):n\in\mathbb N\}}$.}
\end{example} 
\noindent{\bf Acknowledgement.} The authors would like to thank
the referee for his valuable remarks
and comments.

\end{document}